\documentstyle[amssymb]{article}

\newcommand{\proof}{{\bf Proof:  }}
\newcommand{\remark}{{\bf Remark:  }}
\newcommand{\remarks}{{\bf Remarks:  }}
\newcommand{\example}{{\bf Example:  }}
\newcommand{\examples}{{\bf Examples:  }}
\newcommand{\question}{{\bf Question:  }}
\newcommand{\questions}{{\bf Questions:  }}

\newcommand{\charr}{\mbox{char }}

\newcommand{\hb}{\newline\hspace*{\fill}$\Box$}

\newcommand{\SS}{\scriptsize}
\newcommand{\Imm}{\mbox{Im }}
\newcommand{\Ker}{\mbox{Ker }}

\newcommand{\ses}[3]
{\mbox{$0 \rightarrow #1 \rightarrow #2 \rightarrow #3 \rightarrow 0$}}
\newcommand{\codim}{\mbox{codim }}
\newtheorem{theorem}{Theorem}[section]
\newtheorem{lemma}[theorem]{Lemma}
\newtheorem{definition}[theorem]{Definition}
\newtheorem{proposition}[theorem]{Proposition}
\newtheorem{corollary}[theorem]{Corollary}

\begin{document}

\parindent0pt

\title{\bf The monoid of families of quiver representations}
\author{Markus Reineke\\[2ex] BUGH Wuppertal, Gau\ss str. 20, D-42097
Wuppertal, Germany\\ (e-mail: reineke@math.uni-wuppertal.de)}
\date{}
\maketitle

\section{Introduction}

After the work of V. G. Kac (\cite{Ka}), it became clear that the structure
theory of Kac-Moody algebras and the algebro-geometric study of varieties
of representations are valuable tools for the study of representations of 
quivers. This insight was even strengthened by C. M. Ringel's Hall algebra
approach to quantum groups (\cite{Ri2}).\\[1ex]
The aim of this paper is to develop
an algebraic tool, namely the monoid structure on families of quiver representations as stated in the title, which is directly related both to quantized Kac-Moody algebras, and to geometric properties of varieties of representations.
Although this new structure is far from being understood in detail, it it
easily accessible to explicit calculations, and it produces a lot of interesting new examples for natural families of quiver representations.\\[1ex]
The motivation for the definition of this structure comes from the paper 
\cite{Re}, where the author defined a monoid structure on representations of
Dynkin quivers given by so-called generic extensions.
This structure turned out to be interesting for two reasons: on
the one hand, it
provides a computational tool for controlling extensions of two given representations. On the other hand, this monoid of generic extensions relates
to quantum groups: its monoid ring is isomorphic to a degenerate version
of the quantized enveloping algebra (see \cite{KT}) corresponding to the quiver.\\[1ex]
For these reasons, it was desirable to extend this monoid construction
to arbitrary quivers. The main obstacle is
the non-existence of generic extensions for non-Dynkin quivers. Intuitively,
one gets the idea to extend not only individual representations, but whole
families of them to overcome this obstacle.\\[1ex]
In the present paper, this idea is made precise: we define a monoid
structure on irreducible, closed subvarieties of the varieties of representations of an arbitrary quiver without cycles by taking a 'variety
of extensions' as their product. For Dynkin quivers, this notion coincides with the one used in \cite{Re}.\\[1ex]
To retain the relation to quantum groups, it is neccessary to consider the
submonoid generated by simple representations. This is motivated by an
analogous method in Hall algebra theory (see \cite{Ri2}), where the so-called composition
algebra (the subalgebra of the Hall algebra generated by the simple
representations) is isomorpic to a quantized Kac-Moody algebra by \cite{Gr}.\\[1ex]
In fact, the monoid ring of this submonoid turns out to be a factor of a degenerate quantized
Kac-Moody algebra. But in contrast to the Dynkin case, it seems to be a
proper factor in general. This is a rather surprising fact, which remains
to be understood.\\[1ex]
But this submonoid, called the composition monoid, is also
interesting for other reasons. Its elements correspond to families of
quiver representations posessing a composition series
with prescribed factors {\it in prescribed order}.\\[1ex]
Although such families
are a very natural object, it seems that they were never considered before.
The composition monoid thus provides an algebraic tool for their
study: relations in this monoid yield a wealth of information on identities,
inclusions, etc. of these families.\\[2ex]
Due to these results, the original use of the monoid of generic extensions
and its generalizations is complemented by its role in the study of the
geometry of quiver representations. An interesting interaction of algebraic
and geometric methods evolves: on the one hand, the extension monoid and
the composition monoid provide a rich source of examples of
natural families of quiver representations, and they allow some of their
geometric properties to be controlled.\\[1ex]
On the other hand, geometric results, especially the generic properties of representations as studied
in \cite{Sc}, translate into algebraic facts shading light on the
quite mysterious nature of the composition monoid.\\[2ex]
This paper is organized as follows:\\[1ex]
In section 2, we develop some tools which allow us to define the extension monoid (Definition \ref{extmon}). We show that it coincides with the monoid
of generic extensions for Dynkin quivers and give some basic examples. Under the assumption of zero characteristic of the ground field, we develop a formula for the behaviour of the dimension of a 'family' of representations under the product structure (Proposition \ref{mts2}, Corollary 
\ref{cor1}), generalizing a result of A. Schofield from \cite{Sc}. This formula and its
complements (Theorem \ref{cor2}, Proposition \ref{extest}) are the
main tools for the more detailed study in section 5.\\[1ex]
In section 3, we define the composition monoid (Definition \ref{compmon}).
We analyse its elements, the varieties of representations posessing a
composition series of prescribed type, giving a description in terms of
matrices (Lemma \ref{lades}), a dimension estimate (Theorem \ref{dimest}),
and a criterion for inclusions of such families (Lemma \ref{inc}). Various
examples of their behaviour are described.\\[1ex]
The main result of section 4 is the existence of a surjective comparison map (Proposition \ref{rqg}) between a degenerate version of a quantized
Kac-Moody algebra (Lemma \ref{dqkma}) and the monoid ring of the composition monoid. The surprising differences between the Dynkin case and the general
case are discussed (Theorem \ref{rqgd}).\\[1ex]
Section 5 presents a more detailed study of the structure of the composition
monoid. We show how the dimension formulae of section 2 and results
from \cite{Sc} provide us with
a wealth of relations holding in this monoid (Proposition \ref{obs}). Several
explicit types of relations are developed (Corollaries \ref{drel1}, \ref{drel2}, \ref{drel3}). Finally, we develop a 'partial normal form' for elements of
the composition monoid (Theorem \ref{pnf}), using tools from the geometric
study of quiver representations like Schur roots and the generic decomposition
of a dimension vector.\\[2ex]
This work was done while the author  participated in the TMR network `Algebraic Lie Representations' (TMR ERB FMRX-CT97-0100).\\[1ex]
I would like to thank K. Bongartz, C. M. Ringel, A. Schofield and A. Zelevinsky for interesting discussions concerning this paper.

\section{A monoid structure on families of quiver representations}

For all standard notations and results in the representation theory of
quivers and in quantum group theory, we refer the reader to the books \cite{Ri1}, \cite{Lu}, and to \cite{Ka}, \cite{Sc} for geometric aspects.\\[1ex]
Let $Q$ be a finite quiver without oriented cycles. We denote its set of vertices by $I$. The free abelian group ${\bf Z}I$ generated by coordinate vectors $\sigma_i$ for $i\in I$ carries
a non-symmetric inner product, the Euler form of $Q$, which is defined by
$$\langle d, e\rangle:=\sum_{i\in I}d_ie_i-\sum_{i\rightarrow j}d_ie_j$$
for $d,e\in{\bf Z}I$. Here and in the following, the sum $\sum_{i\rightarrow j}$ is meant to run over all arrows in $Q$.\\[1ex]
Let $k$ be an algebraically closed field of characteristic $0$. For $d\in{\bf N}I$, we denote by
$$R_d:=\bigoplus_{\alpha:i\rightarrow j}{\rm Hom}_k(k^{d_i},k^{d_j})$$
the affine $k$-variety of $k$-representations of $Q$ of dimension vector $d$.
In the following, we will always identify $k$-points $M=(M_\alpha)_\alpha$ of $R_d$ with the
corresponding representations of $Q$. The variety $R_d$ is acted upon by the
connected reductive algebraic group 
$$G_d:=\prod_{i\in I}{\rm GL}_{d_i}(k)$$
via
$$(g_i)_i (M_\alpha)_\alpha:=(g_{h\alpha}M_\alpha g_{t\alpha}^{-1})_\alpha,$$
where for an arrow $\alpha$ in $Q$ from $i$ to $j$, we denote by $t\alpha$ and
$h\alpha$ its tail and head, respectively.\\[1ex]
Note that, by definition, we have
$$\langle d,d\rangle=\dim G_d-\dim R_d.$$
The orbits ${\cal O}_M$ under the above action correspond bijectively to
the isomorphism classes $[M]$ of representations of $Q$ of dimension vector $d$.

\begin{definition} For subsets ${\cal A}\subset R_d$, ${\cal B}\subset R_e$, we
define ${\cal A}*{\cal B}\subset R_{d+e}$ to be the subset of all extensions
of representations in ${\cal A}$ by representations in ${\cal B}$, i.e.:
\begin{eqnarray*}
{\cal A}*{\cal B}&=\{&X\in R_{d+e}\, :\, \mbox{ there exists an exact sequence }\\
&&\ses{N}{X}{M} \mbox{ for some } M\in{\cal A}, N\in{\cal B}\}.
\end{eqnarray*}
\end{definition}

\begin{lemma}\label{lem1} If ${\cal A}\subset R_d$, ${\cal B}\subset R_e$ are irreducible,
closed subvarieties, which are stable under the corresponding group actions, then the
same holds for ${\cal A}*{\cal B}$.
\end{lemma}

\proof Denote by $Z_{d,e}$ the closed subvariety of $R_{d+e}$ consisting of
representations of the form
$$(\left[\begin{array}{cc}N_\alpha&\zeta_\alpha\\ 0&M_\alpha\end{array}\right])_\alpha,$$
where $M=(M_\alpha)_\alpha$ and $N=(N_\alpha)_\alpha$ belong to $R_d$ and
$R_e$, respectively, and $\zeta_\alpha$ belongs to $k^{d_{t\alpha}\times
e_{h\alpha}}$ for all arrows $\alpha$. Via the obvious projection $p:Z_{d,e}\rightarrow R_d\times R_e,$
the variety $Z_{d,e}$ becomes a trivial
vector bundle of rank $\sum_{i\rightarrow j}d_ie_j$ over $R_d\times R_e$. Furthermore,
the group action of $G_{d+e}$ on $R_{d+e}$ induces a morphism
$m:G_{d+e}\times Z_{d,e}\rightarrow R_{d+e}.$
Thus,  we have the following diagram of varieties and morphisms:
$$\begin{array}{clcrc}
&&G_{d+e}\times Z_{d,e}&&\\
&\swarrow pr_2&&m\searrow&\\
Z_{d,e}&&&&R_{d+e}.\\
p\downarrow&&&&\\
R_d\times R_e&&&&
\end{array}$$
By definition of an extension between representations, we have:
$${\cal A}*{\cal B}=m(pr_2^{-1}(p^{-1}({\cal A}\times{\cal B}))).$$
Thus, ${\cal A}*{\cal B}$ is nothing else than the $G_{d+e}$-saturation of
$p^{-1}({\cal A}\times{\cal B})$. It follows that ${\cal A}*{\cal B}$ is
always $G_{d+e}$-stable.\\[1ex]
Since $p$ is a trivial bundle, this also proves that ${\cal A}*{\cal B}$
is irreducible, provided ${\cal A}$ and ${\cal B}$ are so. To prove that
${\cal A}*{\cal B}$ is closed, we just have to note that $p^{-1}({\cal A}
\times{\cal B})$ is closed in $R_{d+e}$, and that it is stable under the induced
action of the parabolic
$$P_{d,e}:=\left[\begin{array}{cc}G_e&*\\ 0&G_d\end{array}\right]\subset G_{d+e}.$$ \hb

\begin{lemma} Given subsets ${\cal A}\subset R_d$, ${\cal B}\subset R_e$, ${\cal C}\subset R_f$, we have
$$({\cal A}*{\cal B})*{\cal C}={\cal A}*({\cal B}*{\cal C}).$$ 
\end{lemma}

\proof Both terms in the above equation are easily seen to equal the set of all
representations $X\in R_{d+e+f}$ posessing a filtration $X\supset Y\supset Z$
such that
$$X/Y\in{\cal A},\; Y/Z\in{\cal B},\; Z\in{\cal C}.$$\hb

These two lemmas allow us to make the following definition:

\begin{definition}\label{extmon} The set ${\cal M}(Q)$ of all irreducible, closed, $G_d$-stable
subvarieties of all $R_d$ for $d\in{\bf N}I$, together with the operation $*$
and the unit element $R_0$, is called the extension monoid of $Q$.
\end{definition}

Note that ${\cal M}(Q)$ is naturally ${\bf N}I$-graded by setting
$${\cal M}(Q)_d:=\{{\cal A}\in{\cal M}(Q)\, :\, {\cal A}\subset R_d\}.$$

\examples $\;$
\begin{enumerate}
\item Suppose that $Q$ is of Dynkin type, i.e. the corresponding unoriented
graph is a disjoint union of Dynkin diagrams of type $A_n$, $D_n$, $E_6$,
$E_7$, $E_8$. Then, by Gabriel's Theorem, $G_d$ acts on $R_d$ with a finite
number of orbits for all $d\in{\bf N}I$. Thus, the closed, irreducible, $G_d$-stable subvarieties of the $R_d$ are precisely the closures $\overline{{\cal O}_M}$ of orbits. From \cite{Re}, section 2, it is clear that
$$\overline{{\cal O}_M}*\overline{{\cal O}_M}=\overline{{\cal O}_{M*N}},$$
where $M*N$ denotes the so-called generic extension of $M$ by $N$, i.e. the
unique (up to isomorphism) extension $X$ of $M$ by $N$ with minimal $k$-dimension of its endomorphism ring ${\rm End}_{kQ}(X)$.\\[1ex]
In this case, the monoid ${\cal M}(Q)$ was related in \cite{Re} to degenerate
forms of quantized enveloping algebras. In fact, the definition \ref{extmon}
arose from the wish to generalize the methods of \cite{Re} to arbitrary quivers.
\item Assume that $Q$ is the simplest non-Dynkin quiver, namely the Kronecker
quiver $1{\rightarrow\atop\rightarrow}2$. The varieties $R_{\sigma_1}$ amd
$R_{\sigma_2}$ consist of a single point, and we have
$$R_{\sigma_2}*R_{\sigma_1}=\{0\}\subset R_{\sigma_1+\sigma_2}\mbox{ and }R_{\sigma_1}*R_{\sigma_2}=R_{\sigma_1+
\sigma_2}.$$
But the variety $R_{\sigma_1+\sigma_2}$ is the union of a ${\bf P}^1(k)$-family of
orbit closures intersecting in the single point $0$ corresponding to the
semisimple representation:
$$R_{\sigma_1+\sigma_2}=\bigcup_{\lambda\in{\bf P}^1(k)}\overline{{\cal O}_{M_\lambda}},$$
where in the representation $M_\lambda$, the arrows are represented by
$\lambda_1,\lambda_2$, respectively, for $\lambda=(\lambda_1:\lambda_2)$ in homogeneous coordinates.\\[1ex]
This shows that the concept of generic extensions cannot be generalized to non-Dynkin quivers. Instead, we have to replace individual representations
by `families' of representations, a notion which we formalize by considering
irreducible, closed, $G_d$-stable subvarieties of the $R_d$.
\item Continuing the above example, it is easily seen that
$${\cal M}(Q)_{\sigma_1+\sigma_2}=\{\{0\},R_{\sigma_1+\sigma_2},\overline{{\cal O}_{M_\lambda}}\,
:\, \lambda\in{\bf P}^1(k)\}.$$
The first two elements can be realized as products of other elements of ${\cal M}(Q)$ as shown
above, which is obviously not possible for the $\overline{{\cal O}_{M_\lambda}}$
. This shows that even to write down a set of generators
for ${\cal M}(Q)$ will be difficult in general.
\end{enumerate} 

It is an interesting problem to retrieve geometric information on ${\cal A}*
{\cal B}$ from the geometry of ${\cal A}$ and ${\cal B}$. As a first result
in this direction, we will now develop a formula for the dimension of ${\cal A}
*{\cal B}$. It will be the main tool needed in section 5.\\[1ex]
Given subvarieties ${\cal A}\subset R_d$, ${\cal B}\subset R_e$ from ${\cal M}(Q)$, we
define
$${\rm hom}({\cal B},{\cal A}):=\min\{\dim_k{\rm Hom}_{kQ}(B,A)\, :\,
A\in{\cal A}, B\in{\cal B}\},$$
which is well defined since the function $(A,B)\mapsto\dim{\rm Hom}_{kQ}(B,A)$
is upper semicontinuous on $R_d\times R_e$ by \cite{Bo}. We also define
$${\rm ext}({\cal B},{\cal A}):={\rm hom}({\cal B},{\cal A})-\langle e,d\rangle.$$
To make the formula easier to read, we will formulate it in terms of the
codimension
$$\codim {\cal A}:=\dim R_d-\dim{\cal A}.$$

\begin{proposition}\label{mts2} Given ${\cal A}$ and ${\cal B}$ as above, we have
$$\codim {\cal A}*{\cal B}=\codim {\cal A}+\codim {\cal B}-\langle e,d\rangle+r,$$
where $r\in{\bf N}$ is such that $0\leq r\leq {\rm hom}({\cal B},{\cal A}).$
\end{proposition}

\remark Using the notations to be developed in the proof of this proposition,
we give a precise formula for $r$ in Corollary \ref{cor1}. The formula is
difficult to handle in general; several special cases where $r$ (and ${\rm hom}({\cal B},{\cal A})$) can be computed explicitely are listed in Theorem
\ref{cor2}.\\[2ex]
\proof The main idea of the proof is to generalize the methods of A. Schofield in \cite{Sc},
who considers the special case ${\cal A}=R_d$, ${\cal B}=R_e$. But instead of
passing to a scheme-theoretic viewpoint as in \cite{Sc}, we use the concept of
generic smoothness, as suggested by K. Bongartz.\\[1ex]
Recall the notations of the proof of Lemma \ref{lem1}. We start by choosing subvarieties ${\cal A}'\subset{\cal A}$, ${\cal B}'\subset{\cal B}$ which are smooth, open, and stable under the corresponding
group actions. Consider the subvariety
$$Z'=p^{-1}({\cal A}'\times{\cal B}')\subset Z_{d,e}.$$
By the triviality of $p$, the variety $G_{d+e}\times Z'$ is irreducible, smooth,
and open (thus dense) in $G_{d+e}\times p^{-1}({\cal A}\times{\cal B})$.
So the morphism
$$m:G_{d+e}\times Z'\rightarrow{\cal A}*{\cal B}$$
induced by $m:G_{d+e}\times Z_{d,e}\rightarrow R_{d+e}$ is dominant. By the
theorem on generic smoothness (see \cite{Ha}, III, 10.7), there exists an
open subvariety $V\subset{\cal A}*{\cal B}$, such that
$$m':=m|_{m^{-1}(V)}:m^{-1}(V)\rightarrow V$$
is a smooth morphism. (Note that at this point, we need the standing assumption
$\charr k=0$.) We can also assume this smooth subset $V$ to be
$G_{d+e}$-stable, since $m$ is $G_{d+e}$-equivariant for the actions on $G_{d+e}
\times Z'$ and ${\cal A}*{\cal B}$ given by left multiplication on the left factor and by the natural action, respectively. It follows that $m^{-1}(V)$ is also $G_{d+e}$-stable, hence of the
form $G_{d+e}\times U$ for some open subvariety $U\subset Z'$.\\[1ex]
On the other hand, by the upper semicontinuity of the function $(A,B)\mapsto
\dim{\rm Hom}_{kQ}(B,A)$, we can choose an open subvariety $V'\subset
{\cal A}'*{\cal B}'$ such that $$\dim{\rm Hom}(B,A)={\rm hom}({\cal B},{\cal A})$$
for all $(A,B)\in V'$. Again by the triviality of $p$, the inverse image $U'=p^{-1}(V')$ is open in $Z'$.\\[1ex]
But since $Z'$ is irreducible, the open subsets $U$ and $U'$ have non-zero
intersection; we choose a point
$$z=\left[\begin{array}{cc}B&\zeta\\ 0&A\end{array}\right]\in U\cap U'.$$
By (\cite{Ha}, III, 10.4), the differential
$$dm_{(1,z)}:T_{(1,z)}(G_{d+e}\times U)\rightarrow T_z(V)$$
between the Zariski tangent spaces is surjective since $m'$ is smooth in the point $(1,z)$. We have to compute $dm_{(1,z)}$.\\[1ex]
First, we will view the tangent spaces $T_A{\cal A}$ and $T_B{\cal B}$ as
affine subspaces of the vector spaces $R_d$ and $R_e$, respectively.
The tangent space $T_{(1,z)}(G_{d+e}\times U)$ can thus be viewed as the
affine subspace
$${\frak g}\oplus\left[\begin{array}{cc}T_B{\cal B}&*\\ 0&T_A{\cal A}\end{array}
\right]$$
of ${\frak g}\oplus R_{d+e}$, where ${\frak g}$ denotes the Lie algebra of $G_{d+e}$.\\[1ex]
To compute $dm_{(1,z)}$, we make a calculation using $k[\epsilon]/(\epsilon^2)$-valued
points:
\begin{eqnarray*}
m(1+\epsilon x, z+\epsilon t)&=&(1+\epsilon x)(z+\epsilon t)(1+\epsilon x)^{-1}\\
&=&(1+\epsilon x)(z+\epsilon t)(1-\epsilon x)\\
&=&z+\epsilon (t+xz-zx),
\end{eqnarray*}
thus
$$dm_{(1,z)}(x,t)=t+xz-zx$$
for $x\in{\frak g}$, $t\in\left[\begin{array}{cc}T_B{\cal B}&*\\ 0&T_A{\cal A}\end{array}
\right]\subset R_{d+e}$. (In this and the following formulae, the multiplication of $I$-tuples is component-wise multiplication of matrices).
We compute the kernel of $dm_{(1,z)}$:
\begin{eqnarray*}
\Ker dm_{(1,z)}&=&\{(x,t)\in{\frak g}\oplus\left[\begin{array}{cc}T_B{\cal B}&*\\ 0&T_A{\cal A}\end{array}
\right]\, :\, t=zx-xz\}\\
&\simeq&\{x\in{\frak g}\, :\, zx-xz\in T_zp^{-1}({\cal A}\times{\cal B})\}.
\end{eqnarray*}
Writing
$$x=\left[\begin{array}{ll}x_1&x_2\\ x_3&x_4\end{array}\right]$$
via the identification
$${\frak g}\simeq\prod_{i\in I}\left[\begin{array}{ll}{\rm End}(k^{e_i})&
{\rm Hom}(k^{d_i},k^{e_i})\\
{\rm Hom}(k^{e_i},k^{d_i})&{\rm End}(k^{d_i})\end{array}\right],$$
we can calculate:
\begin{eqnarray*} zx-xz&=&\left[\begin{array}{cc}B&\zeta\\ 0&A\end{array}
\right]\left[\begin{array}{cc}x_1&x_2\\ x_3&x_4\end{array}\right]-\left[\begin{array}{cc}x_1&x_2\\ x_3&x_4\end{array}\right]\left[\begin{array}{cc}B&\zeta\\ 0&A\end{array}
\right]\\
&=&\left[\begin{array}{lc}Bx_1-x_1B+\zeta x_3&*\\ Ax_3-x_3B&Ax_4-x_4A-x_3\zeta
\end{array}\right],
\end{eqnarray*}
thus
\begin{eqnarray*}\Ker dm_{(1,z)}&\simeq&\{\left[\begin{array}{cc}x_1&x_2\\ x_3&x_4\end{array}\right]\in{\frak g}\, :\,Bx_1-x_1B+\zeta x_3\in T_B{\cal B},\\
&&Ax_4-x_4A-x_3\zeta\in T_A{\cal A},\; x_3\in{\rm Hom}_{kQ}(B,A)\}.
\end{eqnarray*}
We see that we have a projection
$$F:\Ker dm_{(1,z)}\rightarrow {\rm Hom}_{kQ}(B,A)$$
given by
$$\left[\begin{array}{cc}x_1&x_2\\ x_3&x_4\end{array}\right]\mapsto x_3.$$
We have
$$\Ker F=\{\left[\begin{array}{cc}x_1&x_2\\ 0&x_4\end{array}\right]\,
:\, Ax_4-x_4A\in T_A{\cal A},\; Bx_1-x_1B\in T_B{\cal B}\}.$$
But the map $x_1\mapsto Bx_1-x_1B$ is just the differential of the map
$$G_d\rightarrow R_d,\;\; g\mapsto gB$$
(analogously for $x_4\mapsto Ax_4-x_4A$), thus the two conditions in the
description of $\Ker F$ are automatically fulfilled.\\[1ex]
This means that $$\Ker F\simeq{\frak p},$$
the Lie algebra of the parabolic
$P_{d,e}\subset G_{d+e}$. Furthermore, we can now identify the image of $F$ as
$$\Imm F=\{f\in{\rm Hom}_{kQ}(B,A)\, :\, -f\zeta\in T_A{\cal A},\;
\zeta f\in T_B{\cal B}\}.$$
We define $r:=\dim\Imm F$. Putting everything together, we can calculate:\\
$\codim{\cal A}*{\cal B}=$
\begin{eqnarray*}
&=&\dim R_{d+e}-\dim{\cal A}*{\cal B}\\
&&\hspace*{10ex}{\mbox{(by definition)}}\\
&=&\dim R_{d+e}-\dim T_z({\cal A}*{\cal B})\\
&&\hspace*{10ex}{\mbox{(since $z\in V$ is a smooth point)}}\\
&=&\dim R_{d+e}-\dim T_{(1,z)}(G_{d+e}\times p^{-1}({\cal A}\times{\cal B}))+
\dim\Ker dm_{(1,z)}\\
&&\hspace*{10ex}{\mbox{(by the surjectivity of $dm_{(1,z)}$)}}\\
&=&\dim R_{d+e}-\dim G_{d+e}-\dim p^{-1}({\cal A}\times{\cal B})+\dim\Ker F+\dim\Imm F\\
&=&\dim R_{d+e}-\dim G_{d+e}-\dim  p^{-1}({\cal A}\times{\cal B})+\dim {\frak p}+r\\
&&\hspace*{10ex}{\mbox{(by the description of $\Ker F$ and $\Imm F$)}}\\
&=&-\langle d+e,d+e\rangle-\dim{\cal A}-\dim{\cal B}-\sum_{i\rightarrow j}d_ie_j+\\
&&+\dim G_d+\dim G_e+\sum_{i\in I}d_ie_i+r\\
&&\hspace*{10ex}{\mbox{(by the properties of $p$)}}\\
&=&-\langle d,d\rangle+\dim G_d-\dim {\cal A}-\langle e,e\rangle+\dim G_e-\dim
{\cal B}-\\
&&-\langle d,e\rangle+\sum_{i\in I}d_ie_i-\sum_{i\rightarrow j}d_ie_j-
\langle e,d\rangle+r\\
&=&\codim {\cal A}+\codim {\cal B}-\langle e,d\rangle+r\\
&&\hspace*{10ex}{\mbox{(using again the
above formula).}}
\end{eqnarray*}
This proves the proposition. \hb

Using the definition of $r$ and following the construction of the special point
$z$ in the proof, we get:

\begin{corollary}\label{cor1} Under the assumptions of Proposition \ref{mts2},
we have
$$r=\dim\{f\in{\rm Hom}_{kQ}(B,A)\, :\, f\zeta\in T_A{\cal A},\; \zeta f\in T_B{\cal B}\}$$
for a tuple $(A,B,\zeta)\in{\cal A}\times{\cal B}\times\bigoplus_{i\rightarrow j}
{\rm Hom}(k^{d_i},k^{e_j})$ such that
\begin{enumerate}
\item $A$ is a smooth point of ${\cal A}$, $B$ is a smooth point of ${\cal B}$, and $\left[\begin{array}{cc}B&\zeta\\ 0&A\end{array}\right]$ is a smooth point of ${\cal A}*{\cal B}$,
\item the morphism $m:G_{d+e}\times p^{-1}({\cal A}\times{\cal B})\rightarrow
{\cal A}*{\cal B}$ is smooth in the point $(1,\left[\begin{array}{cc}B&\zeta\\ 0&A\end{array}\right])$,
\item $\dim{\rm Hom}_{kQ}(B,A)={\rm hom}({\cal B},{\cal A})$.
\end{enumerate}
\end{corollary}

Let us record some special cases where $r$ can be computed.

\begin{theorem}\label{cor2} Under the assumptions of Proposition \ref{mts2},
we have:
\begin{enumerate}
\item $$\codim {\cal A}*{\cal B}\geq \codim {\cal A}+\codim {\cal B}-\langle
e,d\rangle.$$
If ${\rm hom}({\cal B},{\cal A})=0$, equality holds.
\item $$\codim {\cal A}*{\cal B}\leq \codim {\cal A}+\codim {\cal B}+{\rm ext}({\cal B},{\cal A}).$$
If ${\rm ext}({\cal A},{\cal B})=0$, or ${\cal A}=R_d$ and ${\cal B}=R_e$,
equality holds.
\end{enumerate}
\end{theorem}

\proof The inequality of the first part follows immediately from Proposition
\ref{mts2} using $r\geq 0$. If ${\rm hom}({\cal B},{\cal A})=0$, then
${\rm Hom}_{kQ}(B,A)=0$ for a tuple $(A,B,\zeta)$ as in Corollary \ref{cor1}.
So, the explicit formula for $r$ yields $r=0$.\\[1ex]
The inequality of the second part follows immediately from Proposition
\ref{mts2} using $r\leq {\rm hom}({\cal B},{\cal A})$. If ${\rm ext}({\cal A},
{\cal B})=0$, we can assume additionally that ${\rm Ext}^1_{kQ}(A,B)=0$ for
a special tuple $(A,B,\zeta)$ as in Corollary \ref{cor1}. But this means
that we can also assume $\zeta=0$. Thus, the condition on $f$ in the formula
for $r$ becomes void, yielding $r={\rm hom}({\cal B},{\cal A})$.\\[1ex]
The same happens in the case ${\cal A}=R_d$, ${\cal B}=R_e$, since then, the
tangent spaces $T_A{\cal A}$, $T_B{\cal B}$ are just $R_d$, $R_e$, respectively. \hb

\remarks
\begin{enumerate}
\item The second part of Theorem \ref{cor2} immediately implies:
$$R_d*R_e=R_{d+e} \iff {\rm ext}(R_e,R_d)=0,$$
which is Theorem 3.3 of \cite{Sc}.\\[1ex]
This last result was generalized by W. Crawley-Boevey to the case of
arbitrary characteristic of the ground field in \cite{CB}. One can expect
a similar generalization to hold here.
\item It would be interesting to have a similar formula for the 'generic
number of parameters' of ${\cal A}*{\cal B}$, i.e. of
$$\dim {\cal A}*{\cal B}-\max\{\dim {\cal O}_X\, :\, X\in{\cal A}*{\cal B}\}$$
in terms of those for ${\cal A}$ and ${\cal B}$. However, this seems to be
a difficult problem.
\end{enumerate}

Since the dimension formulae make essential use of the values ${\rm ext}({\cal B},{\cal A})$, we add a formula which allows their inductive calculation
in special cases.

\begin{proposition}\label{extest} Given ${\cal A},{\cal B},{\cal C}$ in ${\cal M}(Q)$, we have
$${\rm ext}({\cal A},{\cal B}*{\cal C})\leq{\rm ext}({\cal A},{\cal B})+{\rm ext}({\cal A},{\cal C}).$$
If ${\rm hom}({\cal A},{\cal B})=0$ or ${\rm ext}({\cal B},{\cal C})=0$, then
equality holds.
\end{proposition}

\proof Denote ${\cal D}:={\cal B}*{\cal C}$. We choose representations $A_1,A_2,A_3\in{\cal A}$, $B\in{\cal B}$, $C\in{\cal C}$, $D\in{\cal D}$ such that
$$\dim{\rm Ext}^1(A_1,B)={\rm ext}({\cal A},{\cal B}),\;
\dim{\rm Ext}^1(A_2,C)={\rm ext}({\cal A},{\cal C}),$$
$$\mbox{ and }\dim{\rm Ext}^1(A_3,D)={\rm ext}({\cal A},{\cal D}).$$
Then the set
$${\cal A}_B:=\{A\in{\cal A}\, :\, \dim{\rm Ext}^1(A,B)={\rm ext}({\cal A},{\cal B})\}$$
is open and non-empty in ${\cal A}$. Analogously we define ${\cal A}_C$ and ${\cal A}_D$. Since ${\cal A}$
is irreducible, the intersection ${\cal A}_B\cap{\cal A}_c\cap{\cal A}_D$ is
again non-empty. We choose a representation $A$ from this intersection.
Now we consider the sets
$${\cal B}_A:=\{B\in{\cal B}\, :\, \dim{\rm Ext}^1(A,B)={\rm ext}({\cal A},{\cal B})\}\subset{\cal B},$$
analogously for ${\cal C}_A\subset{\cal C}$ and ${\cal D}_A\subset{\cal D}$.
All these subsets are open and non-empty. Thus, ${\cal B}_A*{\cal C}_A$ is
dense in ${\cal D}$. We choose a representation $X$ from the intersection
$({\cal B}_A*{\cal C}_A)\cap{\cal D}_A$. Summing up this construction, we find representations $A,B,C,X$ such that
$$\dim{\rm Ext}^1(A,B)={\rm ext}({\cal A},{\cal B}),\;
\dim{\rm Ext}^1(A,C)={\rm ext}({\cal A},{\cal C}),$$
$$\mbox{ and }\dim{\rm Ext}^1(A,X)={\rm ext}({\cal A},{\cal D}),$$
and such that there exists an exact sequence
$$\ses{C}{X}{B}.$$
We thus have an induced exact sequence
$${\rm Hom}(A,B)\stackrel{\partial}{\rightarrow}{\rm Ext}^1(A,C)\rightarrow
{\rm Ext}^1(A,X)\rightarrow{\rm Ext}^1(A,B)\rightarrow 0,$$
since the category of representations of a quiver has global dimension at most one. From this sequence, we get the desired estimate
\begin{eqnarray*}
{\rm ext}({\cal A},{\cal B}*{\cal C})&=&\dim{\rm Ext}^1(A,X)\\
&\leq&\dim{\rm Ext}^1(A,B)+\dim{\rm Ext}^1(A,C)\\
&=&{\rm ext}({\cal A},{\cal B})+{\rm ext}({\cal A},{\cal C}).
\end{eqnarray*}
Now assume that ${\rm hom}({\cal A},{\cal B})=0$ (resp. ${\rm ext}({\cal B},{\cal C})=0$). Then we can assume additionally that ${\rm Hom}(A,B)=0$
(resp. ${\rm Ext}^1(B,C)=0$). In both cases, we can assume $\partial$ to
be the zero map, thus equality holds in the above estimate. \hb

\section{The composition monoid}

As the last example following Definition \ref{extmon} already indicates, the
extension monoid ${\cal M}(Q)$ seems to be 'too big' to allow for a
reasonable understanding of its algebraic structure. Therefore, we will now
consider a reasonable submonoid, namely the one generated by the
varieties $R_{\sigma_i}$ for $i\in I$. Note that $R_{\sigma_i}$ is just the orbit of
the simple representation $E_i$ of $Q$ corresponding to the vertex $i$.\\[1ex]
This restriction is also motivated by the analogy of ${\cal M}(Q)$ to
quantized enveloping algebras and Hall algebras (which will be explained in
detail in section 4). In this context, one usually considers the subalgebra of the Hall algebra (the composition algebra, see \cite{Ri2}) which is generated
by the isomorphism classes of simple representations. It is isomorphic to
the quantized enveloping algebra of the Kac-Moody algebra corresponding to $Q$
by \cite{Gr}, whereas the Hall algebra itself is the quantized enveloping algebra of some rather mysterious Borcherds algebra by \cite{SV}.

\begin{definition}\label{compmon} The submonoid of ${\cal M}(Q)$ generated by the $R_{\sigma_i}$
for $i\in I$ is called the composition monoid of $Q$ and is denoted by
${\cal C}(Q)$.
\end{definition}

Obviously, ${\cal C}(Q)$ inherits the ${\bf N}I$-grading from ${\cal M}(Q)$.\\[1ex]
Let $\Omega$ be the set of finite words in the alphabet $I$, which is
${\bf N}I$-graded by defining the degree of the word $(i)$ of length one as $\sigma_i$ for $i\in I$. We denote the degree of an arbitrary word $\omega\in\Omega$ by $|\omega|\in{\bf N}I$.\\[1ex]
Trivially,
we have a surjective map of (graded) monoids
$$\pi:\Omega\rightarrow{\cal C}(Q),\;\; \pi(i)=R_{\sigma_i}.$$
But this map is never an isomorphism. This will be proved in the next section;
for the moment, let us only write down some simple relations between the
$R_{\sigma_i}$ which can be verified directly.

\begin{lemma}\label{rel1} The following relations hold in ${\cal C}(Q)$:
\begin{enumerate}
\item If there is no arrow between $i$ and $j$, then
$$R_{\sigma_i}*R_{\sigma_j}=R_{\sigma_i+\sigma_j}=R_{\sigma_j}*R_{\sigma_i}.$$
\item If there is exactly one arrow from $i$ to $j$, then
$$R_{\sigma_i}*R_{\sigma_j}*R_{\sigma_i}=R_{\sigma_i}*R_{\sigma_i}*R_{\sigma_j}\mbox{ and}$$
$$R_{\sigma_j}*R_{\sigma_i}*R_{\sigma_j}=R_{\sigma_i}*R_{\sigma_j}*R_{\sigma_j}.$$
\end{enumerate}
\end{lemma}

Thus, it becomes neccessary to understand the kernel of the map $\pi$.
Some results in this direction will be the content of section 5.\\[1ex]
We denote the image of a word $\omega\in\Omega$ under the map $\pi$ by
${\cal E}_\omega$. We will now discuss the structure of the subvarieties
${\cal E}_\omega\subset R_{|\omega|}$.\\[1ex]
From the definition of the multiplication $*$ and the fact that $R_{\sigma_i}$ consists of the single orbit ${\cal O}_{E_i}$, the next result follows
immediately.

\begin{lemma} For a word $\omega=(i_1\ldots i_s)\in\Omega$, the irreducible,
closed, $G_{|\omega|}$-stable subvariety ${\cal E}_{\omega}$ of $R_{|\omega|}$
consists of those representations $X$ of $Q$ which admit a composition series
of type $\omega$, i.e. a filtration
$$X=X_0\supset X_1\supset\ldots\supset X_s=0$$
such that $X_{k-1}/X_k\simeq E_{i_k}$ for all $k=1\ldots s$.
\end{lemma}

\remark This description makes it even more important to know the kernel
of $\pi$, since we now know that this provides interesting information
on the existence of composition series of prescribed type.\\[2ex]
We can also describe the varieties ${\cal E}_\omega$ in terms of matrices,
as follows:\\[1ex]
Given a word $\omega=(i_1\ldots i_s)$, we define a function $v:\{1,
\ldots,s\}\rightarrow{\bf N}$ as follows:
\begin{eqnarray*} v(k)=1&\mbox{ if }& i_l\not= i_k\mbox{ for all }l>k\\
v(k)=v(l)+1&\mbox{ if }&l>k,\; i_l=i_k,\; i_p\not= i_k\mbox{ for all }
k<p<l.
\end{eqnarray*}
In other words: each $i$ appears at position $k$ for the $v(k)$-th time,
reading $\omega$ from right to left.
\begin{lemma}\label{lades} The variety ${\cal E}_\omega$ consists of all
representations of $Q$
which are conjugate to a representation $X$ satisfying the following:\\
for all $\alpha\in Q_1$ and $1\leq k<l\leq s$ such that $i_k=t(\alpha)$, $i_l=s(\alpha)$, we have
$$(X_\alpha)_{v(k),v(l)}=0.$$
\end{lemma}

\proof Recall from Lemma \ref{lem1} that ${\cal A}*{\cal B}$ is the saturation
of $p^{-1}({\cal A}\times{\cal B})$ under the natural group action. Applying
this fact repeatedly, we see that ${\cal E}_\omega=R_{\sigma_{i_1}}*\ldots
*R_{\sigma_{i_s}}$ is the $G_{|\omega|}$-saturation of the set of all representations
$X$ in $R_{|\omega|}$ of the form
$$X=\left[\begin{array}{cccc}(E_{i_s})_\alpha&&&*\\
&(E_{i_{s-1}})_\alpha&&\\
&&\ldots&\\
0&&&(E_{i_1})_\alpha
\end{array}\right].$$
Now we note that the formats of the matrices representing $E_i$ are given by
$$(E_i)_\alpha\in\left\{\begin{array}{lcl}
k^{0\times 0}&,&t\alpha\not= i\not= h\alpha,\\
k^{0\times 1}&,&t\alpha= i\not= h\alpha,\\
k^{1\times 0}&,&t\alpha\not= i= h\alpha.
\end{array}\right.$$
Simplifying the above shape of $X$ using this information, we arrive easily
at the description of $X$ as in the statement of the lemma. \hb

\example Consider the quiver
$$i\stackrel{\alpha}{\rightarrow}j\stackrel{\beta}{\leftarrow}k$$
and the word
$$\omega=(iikijkkjjikijijjjkkij).$$
Then we have
$$v=(766585476433524322111),$$
and thus ${\cal E}_\omega$ consists of all representations conjugate to
representations $X$ of the form
$$\SS
X_\alpha=\left[\begin{array}{ccccccc}
*&*&*&*&*&*&*\\
0&*&*&*&*&*&*\\
0&*&*&*&*&*&*\\
0&*&*&*&*&*&*\\
0&0&*&*&*&*&*\\
0&0&0&0&*&*&*\\
0&0&0&0&*&*&*\\
0&0&0&0&*&*&*
\end{array}\right],\;\;\;
X_\beta=\left[\begin{array}{cccccc}
*&*&*&*&*&*\\
0&0&*&*&*&*\\
0&0&*&*&*&*\\
0&0&*&*&*&*\\
0&0&*&*&*&*\\
0&0&0&*&*&*\\
0&0&0&*&*&*\\
0&0&0&0&0&*
\end{array}\right].$$

From the above description of the subvarieties ${\cal E}_\omega$, it is clear
that they are a very natural object to study when considering the geometry of
quiver representations. As a first step, we give an estimate for their
(co-)dimension using Theorem \ref{cor2}.

\begin{theorem}\label{dimest} Assume $\omega=(i_1\ldots i_s)$ is a word in $\Omega$. Then
$$\codim {\cal E}_\omega\geq-\sum_{i\in I}|\omega|_i(|\omega|_i-1)/2+
\sum_{1\leq k<l\leq s}\#\{\mbox{arrows from $i_l$ to $i_k$}\}.$$
\end{theorem}

\proof Repeatedly applying the inequality of the first part of
Theorem \ref{cor2}, we get:
\begin{eqnarray*}
\codim{\cal E}_\omega&=&\codim R_{\sigma_{i_1}}*\ldots*R_{\sigma_{i_s}}\\
&\geq&\sum_{k=1}^s\underbrace{\codim R_{\sigma_{i_k}}}_{=0}-\sum_{1\leq k<l
\leq s}\langle \sigma_{i_l},\sigma_{i_k}\rangle\\
&=&-\sum_{{1\leq k<l\leq s}\atop{i_k=i_l}}\langle \sigma_{i_l},\sigma_{i_k}\rangle-
\sum_{{1\leq k<l\leq s}\atop{i_k\not=i_l}}\langle \sigma_{i_l},\sigma_{i_k}\rangle\\
&=&-\#\{(k,l)\, :\, 1\leq k<l\leq s,\; i_k=i_l\}-
\sum_{1\leq k<l\leq s}\langle \sigma_{i_l},\sigma_{i_k}\rangle\\
&=&-\sum_{i\in I}|\omega|_i(|\omega|_i-1)/2+
\sum_{1\leq k<l\leq s}\#\{\mbox{arrows from $i_l$ to $i_k$}\},
\end{eqnarray*}
where we could drop the assumption $i_k\not= i_l$ from the second sum since
there are no loops in $Q$. \hb

\remark Surprisingly, this simple estimate is in fact an equality in a lot of examples. It would be important to know a sufficiently general class of words where 
equality could be proved.\\[2ex]
Now we discuss inclusions and intersections of the varieties ${\cal E}_\omega$.\\[1ex]
We choose a total ordering on $I$ such that the existence of an arrow from $i$ to $j$ implies $i<j$ (note that this is possible since $Q$ has no oriented cycles). This order is extended to a partial ordering on each $\Omega_d$ for $d\in{\bf N}I$ by the following definition.

\begin{definition} Given $\omega,\omega'\in\Omega_d$, we say that $\omega\leq \omega'$ if there exists a sequence of words
$$\omega=\omega_0,\omega_1,\ldots,\omega_t=\omega'$$
such that for all $k=1\ldots t$, there exist words $\omega_l,\omega_r\in\Omega$ and vertices $i<j$ in $I$ such that
$$\omega_{k-1}=\omega_l\, (ij)\,\omega_r,\;\; \omega_k=\omega_l\, (ji)\,\omega_r.$$
\end{definition}

With this definition, we have:

\begin{lemma}\label{inc} If $\omega\leq\omega'$, then ${\cal E}_\omega\supset{\cal E}_{\omega'}$.
\end{lemma}

\proof Choose a sequence of words $\omega=\omega_0,\omega_1,\ldots,\omega_t=\omega'$ as in the definition.
Then, for all $k=1\ldots t$, we have:
\begin{eqnarray*}
{\cal E}_{\omega_{k-1}}&=&{\cal E}_{\omega_l}*R_{\sigma_i}*R_{\sigma_j}*{\cal E}_{\omega_r}\\
&\supset&{\cal E}_{\omega_l}*\{0\}*{\cal E}_{\omega_r}\\
&=&{\cal E}_{\omega_l}*R_{\sigma_j}*R_{\sigma_i}*{\cal E}_{\omega_r}\\
&=&{\cal E}_{\omega_k},
\end{eqnarray*}
where the inclusion is obvious: since $i<j$, there is no arrow
from $j$ to $i$, so $R_{\sigma_j}*R_{\sigma_i}=\{0\}\subset R_{\sigma_i+\sigma_j}$. \hb

The converse statement does not hold, as can already be seen in
simple examples.\\[1ex]
\examples
\begin{enumerate}
\item If $i<j$, and there is no arrow between $i$ and $j$, then ${\cal E}_{(ji)}={\cal E}_{(ij)}$, but $(ji)\not\leq (ij)$.
\item A more interesting example is the following: suppose there is a single
arrow from $i$ to $j$. Then, using Lemma \ref{rel1}, we have
$${\cal E}_{(jiiij)}={\cal E}_{(jiiji)}\subset{\cal E}_{(ijiji)}$$
since $(ijiji)<(jiiji)$, but $(ijiji)\not\leq (jiiij)$.
\end{enumerate}

Nevertheless, we can ask whether the following partial converse to Lemma \ref{inc} holds:\\[1ex]
\question Suppose ${\cal E}_\omega\supset {\cal E}_{\omega'}$. Does it then
follow that there exist words $\overline{\omega},\overline{\omega}'$ such that
$${\cal E}_\omega={\cal E}_{\overline{\omega}},\;{\cal E}_{\omega'}={\cal E}_{\overline{\omega}'}\mbox{ and }\overline{\omega}\leq\overline{\omega}'?$$\\[1ex]
Turning to intersections of the ${\cal E}_\omega$, we construct an example
where ${\cal E}_\omega\cap{\cal E}_{\omega'}$ is not of the form ${\cal E}_{\omega''}$.\\[1ex]
\example Consider the quiver
$$Q\; :\; i{\rightarrow \atop\rightarrow}j{\rightarrow
\atop\rightarrow}k.$$
We order the vertices by $i<j<k$ and consider the dimension
vector $d=\sigma_i+2\sigma_j+\sigma_k$. The Hasse diagram of the poset $\Omega_d$ is
$$\begin{array}{ccccccccc}&&&&(ijjk)&&&&\\
&&&\swarrow&&\searrow&&&\\
&&(jijk)&&&&(ijkj)&&\\
&\swarrow&&\searrow&&\swarrow&&\searrow&\\
(jjik)&&&&(jikj)&&&&(ikjj)\\
\downarrow&&&&\downarrow&&&&\downarrow\\
(jjki)&&&&(jkij)&&&&(kijj)\\
&\searrow&&\swarrow&&\searrow&&\swarrow&\\
&&(jkji)&&&&(kjij)&&\\
&&&\searrow&&\swarrow&&&\\
&&&&(kjji).&&&&
\end{array}$$
A representation of $Q$ of dimension vector $d$ is given by a tuple
$(v_1,v_2,\phi_1,\phi_2)$, where $v_1,v_2\in k^{2\times 1}$, $\phi_1,
\phi_2\in k^{1\times 2}$. We represent such a tuple as a pair of $2\times 2$-matrices
$$A=\left[v_1 v_2\right],\;\; B=\left[{\phi_1\atop\phi_2}\right],$$
on which the group $G_d=k^*\times{\rm Gl}_2(k)\times k^*$ acts via
$$(\lambda,g,\mu)(A,B)=(\frac{1}{\lambda}gA,\mu Bg^{-1}).$$
Using e.g. Lemma \ref{lades}, it is easily verified that
$${\cal E}_{(jijk)}=\{(A,B)\, :\, \det A=0\},\;\;{\cal E}_{(ijkj)}=\{(A,B)\, :\, \det B=0\}$$
and
$${\cal E}_{(jikj)}=\{(A,B)\, :\, \det A=0=\det B,\; BA=0\}.$$
Thus, we have a strict inclusion ${\cal E}_{(jijk)}\cap{\cal E}_{(ijkj)}\supset 
{\cal E}_{(jikj)}$, finishing the example.\\[3ex]
The ascending chain of closed subvarieties ${\cal E}_{\omega}\, :\, \omega
\in\Omega_d$ of $R_d$ makes it possible to single out a 'canonical' open
subset of $R_d$:

\begin{definition} For each $d\in{\bf N}I$, the open subvariety $S_d$ of $R_d$
is defined as the complement of all ${\cal E}_{\omega}\not= R_d$.
\end{definition}

We can ask, for example, whether the categorial quotient $S_d/G_d$ (in the sense
of Mumford) is an interesting object, or under what circumstances a geometric
quotient exists.\\[1ex]
\examples
\begin{enumerate}
\item If $Q$ is a Dynkin quiver, then we will see in the next section that
$S_d$ is precisely the open orbit, so we always have a geometric quotient
consisting of a single point.
\item If $Q$ is the Kronecker quiver, then $S_{\sigma_1+\sigma_2}/G_{\sigma_1+\sigma_2}\simeq{\bf P}^1(k)$ and $$S_{2\sigma_1+2\sigma_2}/G_{2\sigma_1+2\sigma_2}\simeq{\bf P}^1(k)\times{\bf P}^1(k),$$
but the latter is not a geometric quotient, since the fibre of the quotient
map above the diagonal consists of two orbits, corresponding to a four-
dimensional indecomposable representation and a direct sum of two copies of
a two-dimensional indecomposable representation.
\item If $Q$ is the quiver $i{\rightarrow \atop\rightarrow}j{\rightarrow
\atop\rightarrow}k$, then, in the notations of the example above, we have
$$S_d=\{(A,B)\in (k^{2\times 2})^2\, :\, \det A\not=0\not=\det B\}.$$
An easy calculation shows that $S_d/G_d$ is isomorphic to ${\bf P}^3(k)\setminus Q$, where $Q=\{(a:b:c:d)\in{\bf P}^3(k)\, :\, ad=bc\}$.
\end{enumerate}

\section{Relation to quantum groups}

In this section, we show how ${\cal C}(Q)$ relates to quantized enveloping
algebras and Hall algebras.\\[1ex]
Given nonnegative integers $M$ and $N$, we introduce the $q$-factorial and 
$q$-binomial coefficients
$$[M]!=(q-1)\cdot\ldots\cdot(q^{M}-1),\;\;\;
[{{M+N}\atop{M}}]=\frac{[M+N]!}{[M]![N]!}.$$

\begin{definition} Let ${\cal U}^+(Q)$ be the ${\bf Q}[q]$-algebra with
generators $E_i$ for $i\in I$ and relations
$$\sum_{p+p'=n+1}(-1)^{p'}q^{p'(p'-1)}\left[{{p+p'}\atop{p}}\right]E_i^pE_jE_i^{p'}=0,$$
$$\sum_{p+p'=n+1}(-1)^{p}q^{p(p-1)}\left[{{p+p'}\atop{p}}\right]E_j^pE_iE_j^{p'}=0,$$
if there is no arrow from $j$ to $i$, and $n$ is the number of arrows from
$i$ to $j$.
\end{definition}

This algebra is ${\bf N}I$-graded by defining the degree $|E_i|$ of the
generator $E_i$ as $\sigma_i\in{\bf N}I$. This is possible, since the defining
relations above are obviously homogeneous for this grading.\\[1ex]
We consider three variants of this algebra:\\[1ex]
First, we view the function field ${\bf Q}(v)$ as a ${\bf Q}[q]$-algebra via $q=v^2$. Then, we define ${\cal U}^+_v(Q)$ as the scalar extension
${\bf Q}(v)\otimes_{{\bf Q}[q]}{\cal U}^+(Q)$, together with the
twisted multiplication
$$x*y:=v^{\langle |x|,|y|\rangle}x\cdot y$$
for homogeneous elements $x,y$. If we apply this twist to the defining
relations of ${\cal U}^+(Q)$, we see that ${\cal U}^+_v(Q)$ is defined
by the relations
$$\sum_{p+p'=n+1}(-1)^{p'}E_i^{(p)}E_jE_i^{(p')}=0,$$
where $n$ denotes the number of arrows between $i$ and $j$, and $E_i^{(n)}
$ denotes the usual quantized divided powers. By (\cite{Lu}, 33.1.4), these
are precisely the defining relations of the quantized enveloping algebra of the Kac-Moody algebra of the unoriented diagram corresponding to $Q$.\\[1ex]
Next, we can specialize the variable $q$ to any prime power to obtain
${\bf Q}$-algebras ${\cal U}^+_q(Q)$. Combining \cite{Ri2} and \cite{Gr}, these
algebras are precisely the Hall algebras $H({\bf F}_qQ)$ of the quiver $Q$ over
the finite fields ${\bf F}_q$ with $q$ elements.\\[1ex]
Finally, we can specialize the variable $q$ to $0$ to obtain the ${\bf Q}$-algebra
${\cal U}^+_0(Q)$. 

\begin{lemma}\label{dqkma} ${\cal U}^+_0(Q)$ is defined by generators $E_i$ for $i\in I$ and relations
$$E_i^{n+1}E_j=E_i^nE_jE_i,\;\;\; E_iE_j^{n+1}=E_jE_iE_j^n,$$
if there is no arrow from $j$ to $i$, and $n$ is the number of arrows from
$i$ to $j$.
\end{lemma}

\proof We have to consider the constant term of the polynomial $$(-1)^{p'}q^{p'(p'-1)}\left[{{p+p'}\atop{p}}\right]\mbox{ (resp. } (-1)^{p}q^{p(p-1)}\left[{{p+p'}\atop{p}}\right]\mbox{)}$$
appearing in the defining
relations for ${\cal U}^+(Q)$. By definition, the constant term of $
\left[{{p+p'}\atop{p}}\right]$ always equals $1$. Thus, the only non-zero
constant terms appear for $p'$ (resp. $p$) equal to $0$ or $1$. In this case
we arrive immediately at the claimed relations in ${\cal U}^+_0(Q)$. \hb

Borrowing a result from the next section, we get
the promised relation between quantized enveloping algebras and the
composition monoid:

\begin{proposition}\label{rqg} The assignment $E_i\mapsto R_{\sigma_i}$ for $i\in I$
extends uniquely to a surjective homomorphism of ${\bf Q}$-algebras
$$\eta:{\cal U}^+_0(Q)\rightarrow {\bf Q}{\cal C}(Q),$$
where ${\bf Q}{\cal C}(Q)$ denotes the monoid ring of the monoid ${\cal C}(Q)$.
\end{proposition}

\proof By Corollary \ref{drel1}, the defining relations of ${\cal U}^+_0(Q)$
also hold in ${\cal C}(Q)$. Thus, $\eta$ is a well-defined homomorphism of ${\bf Q}$-algebras. It is surjective since 
${\cal C}(Q)$ is generated by the $R_{\sigma_i}$ for $i\in I$ by definition.  \hb

Once this is known, it is reasonable trying to understand ${\cal C}(Q)$ by
studying the comparison map $\eta$, since this possibly allows for
application of quantum group techniques. A basic result in this direction
is the following:

\begin{theorem}\label{rqgd} If $Q$ is of Dynkin type, then the map $\eta$ is an isomorphism. Moreover, we have ${\cal C}(Q)={\cal M}(Q)$ in this case.
\end{theorem}

\proof Assume that $Q$ is of Dynkin type. From the first example following Definition \ref{extmon}, we know that
${\cal M}(Q)$ is isomorphic to the monoid of generic extensions of \cite{Re}.
The composition of the comparison map $\eta$ and the natural inclusion ${\bf Q}{\cal C}(Q)\subset{\bf Q}{\cal M}(Q)$ maps the generator $E_i$ to the
orbit of the simple representation $E_i$. But this map is shown to be an
algebra isomorphism in (\cite {Re}, Theorem 4.2). Thus, both $\eta$ and
the inclusion are actually isomorphisms themselves. \hb

It is reasonable to assume that Dynkin quivers are in fact the only
quivers where $\eta$ is an isomorphism. This is motivated by the
following example:\\[2ex]
\example Let $Q$ be the Kronecker quiver. By Corollary \ref{drel2}, we have
$$R_{\sigma_1+\sigma_2}^{*2}=R_{2\sigma_1+2\sigma_2}.$$
But in ${\cal U}^+_0(Q)$, we have $$E_1E_2E_1E_2\not=E_1^2E_2^2,$$
since the defining relations are of degree $3\sigma_1+\sigma_2$ and $\sigma_1+3\sigma_2$, respectively.\\[2ex]
In fact, one can expect this example to generalize to arbitrary extended
Dynkin quivers. But then, it already generalizes to arbitrary non-Dynkin quivers, since we can always find a subquiver of extended Dynkin type. 

\section{The structure of the monoid ${\cal C}(Q)$}

In this section, we take a closer look at the relations holding in
${\cal C}(Q)$. We will mainly use Theorem \ref{cor2} and the notions
${\rm hom}(R_d,R_e)$, ${\rm ext}(R_d,R_e)$.\\[1ex]
Recall the total ordering on $I$ from section 3. Without loss of generality,
we can thus assume $I=\{1,\ldots,n\}$ with the natural ordering. Using this
notation, it is easy to see that all varieties $R_d$ belong to ${\cal C}(Q)$:

\begin{lemma} For all $d\in{\bf N}I$, we have
$$R_d=R_{\sigma_1}^{*d_1}*\ldots*R_{\sigma_n}^{*d_n}.$$
In particular, $R_d\in{\cal C}(Q)$.\end{lemma}

Once this is known, the first remark following Theorem \ref{cor2} becomes
important for our study since it produces a wealth of relations in ${\cal C}(Q)$. For this reason, we record this observation as a separate proposition.

\begin{proposition}\label{obs} The following relations hold in ${\cal C}(Q)$:
$$\mbox{If }{\rm ext}(e,d):={\rm ext}(R_e,R_d)=0\mbox{, then }R_d*R_e=R_{d+e}.$$
\end{proposition}

To use this proposition, it is neccessary to have a precise criterion for
the vanishing of ${\rm ext}(e,d)$. Fortunately, such a criterion is provided
by A. Schofield's work, although it is difficult to handle in practice due
to its recursive nature:

\begin{theorem}\label{csc} (\cite{Sc}, Theorem 5.4) For $d,e\in{\bf N}I$, the following statements are equivalent:
\begin{enumerate}
\item ${\rm ext}(e,d)=0$,
\item $\langle e',d\rangle\geq 0$ for all $e'\leq e$ such that $R_{e-e'}*R_{e'}=R_e$,
\item $\langle e,d'\rangle\geq 0$ for all
$d'\leq d$ such that $R_{d'}*R_{d-d'}=R_d$.
\end{enumerate}
\end{theorem}

Nevertheless, we can apply this result in several special cases.
We note two such cases which were used in the previous section to 
establish the relation of ${\cal C}(Q)$ to quantum groups.

\begin{corollary}\label{drel1} The following relations hold in ${\cal C}(Q)$:
$$R_{\sigma_i}^{*(n+1)}*R_{\sigma_j}=R_{\sigma_i}^{*n}*R_{\sigma_j}*R_{\sigma_i},\;\;\; R_{\sigma_i}*R_{\sigma_j}^{*(n+1)}=R_{\sigma_j}*R_{\sigma_i}*R_{\sigma_j}^{*n},$$
if there is no arrow from $j$ to $i$, and $n$ is the number of arrows from
$i$ to $j$.
\end{corollary}

\proof To prove the first relation, let $i,j\in I$ be as above and define $d=n\sigma_i+\sigma_j$, $e=\sigma_i$. We apply the second condition of Theorem \ref{csc}.
The only vectors $e'$ to be considered are obviously $e'=0$ and $e'=\sigma_i$.
Thus, ${\rm ext}(\sigma_i,d)=0$ if and only if $\langle \sigma_i,d\rangle\geq 0$. In
our case, $\langle \sigma_i,d\rangle=0$, so we are done by Proposition \ref{obs}. The second relation is
proved similar, this time using the third condition of Theorem \ref{csc}. \hb

We call a dimension vector $d$ an isotropic root if $d$ belongs to
the positive part of the (possibly infinite) root system of the Kac-Moody algebra associated to $Q$, and $\langle d,d\rangle=0$. Note in particular,
that the imaginary roots of an extended Dynkin diagram are isotropic.

\begin{corollary}\label{drel2} If $d$ is an isotropic root, then for all
$m,n\in{\bf N}$, we have
$$R_{md}*R_{nd}=R_{(m+n)d}.$$
\end{corollary}

\proof By (\cite{Sc}, Theorem 3.6), we have ${\rm ext}(d,d)=0$. Thus, there
exist representations $X,Y\in R_d$ such that ${\rm Ext}^1(X,Y)=0$. But this
implies $${\rm Ext}^1(X^n,Y^m)=0\mbox{, hence }{\rm ext}(nd,md)=0$$ for all $m,n\in{\bf N}$. By Proposition \ref{obs}, we are done. \hb

The role of Proposition \ref{obs} is quite ambivalent: on the one hand, we
have potentially a big number of relations in ${\cal C}(Q)$. On the other hand,
this abundance is difficult to control: several examples show that most of the relations of the given kind are redundant, i.e. they already follow from
relations of smaller degree as in the next example.\\[2ex]
\example Let $Q$ be the quiver with two vertices $i<j$ and three
arrows from $i$ to $j$. Using Theorem \ref{csc}, one can see by a direct calculation that ${\rm ext}(2\sigma_i+3\sigma_j,3\sigma_i+\sigma_j)=0$. But the corresponding
relation $R_{3\sigma_i+\sigma_j}*R_{2\sigma_i+3\sigma_j}=R_{5\sigma_i+4\sigma_j}$ can already be seen
using only the relations of Corollary \ref{drel1}:
\begin{eqnarray*}R_{3\sigma_i+\sigma_j}*R_{2\sigma_i+3\sigma_j}&=&R_{\sigma_i}^{*3}*R_{\sigma_j}*R_{\sigma_i}*R_{\sigma_i+3\sigma_j}\\
&=&R_{\sigma_i}^{*4}*R_{\sigma_j}*R_{\sigma_i+3\sigma_j}\\
&=&R_{4\sigma_i}*R_{\sigma_j}*R_{\sigma_i}*R_{\sigma_j}^{*3}\\
&=&R_{4\sigma_i}*R_{\sigma_i}*R_{\sigma_j}^{*4}\\
&=&R_{5\sigma_i+4\sigma_j}.
\end{eqnarray*}\\[1ex]
To see how Proposition \ref{obs} and Theorem \ref{csc} can be used in
calculations, let us analyse the case $d=\sigma_i+k\sigma_j$ for $1\leq k\leq n$ and
$e=x\sigma_i+y\sigma_j$, where $n$ is the number of arrows from $i$ to $j$, assuming there
is no arrow from $j$ to $i$. Using both the proposition and the theorem, we see
that
$$R_d*R_e=R_{d+e} \iff \langle e,d'\rangle\geq 0$$
for all $d'\leq d$ such that $R_{d'}*R_{d-d'}=R_d$, i.e. ${\rm ext}(d-d',d')=0$,
using the proposition again. We have to distinguish
two cases, namely $d'$ being of the form $l\sigma_j$ or $\sigma_i+l\sigma_j$ for $0\leq l\leq k$. In the first case, the theorem yields
$${\rm ext}(d-d',d')=0 \iff \langle \sigma_i+(k-l)\sigma_j,p\sigma_j\rangle\geq 0
\mbox{ for all }p\leq l,$$
which is only fulfilled in case $l=0$ due to the assumption $k\leq n$. For the second case,
we get
$${\rm ext}(d-d',d')=0 \iff \langle p\sigma_j, \sigma_i+l\sigma_j\rangle\geq 0\mbox{ for all }
p\leq k-l,$$
which is fulfilled for all $l=0\ldots k$. Thus, we have
$$R_d*R_e=R_{d+e} \iff \langle e,d\sigma_i+l\sigma_j\rangle\geq 0\mbox{ for }l=0\ldots k.$$
This last condition is computed to
$$x\geq (nx-y)l\mbox{ for all }l=0\ldots k.$$
We arrive at two cases: if $nx\leq y$, the above condition is obviously fulfilled, thus we get the relations
$$R_{\sigma_i}*R_{\sigma_j}^{*k}*R_{\sigma_i}^{*x}*R_{\sigma_j}^{*y}=R_{\sigma_i}^{*(x+1)}*R_{\sigma_j}^{*(y+k)}\mbox{ if } y\leq nx;$$
but all these relations are easily seen to follow from
$$R_{\sigma_j}*R_{\sigma_i}*R_{\sigma_j}^{*n}=R_{\sigma_i}*R_{\sigma_j}^{*(n+1)}.$$
The other case is more interesting: if $nx>y$, we have the extra condition
$x\geq(nx-y)k$. These inequalities are easily reformulated as 
$$(n-\frac{1}{k})x\leq y<nx.$$
Such an $y$ can only be found in case $x\geq k$, so we get the following new
relations:

\begin{corollary}\label{drel3} Given $1\leq k\leq n$, $x\geq k$, define
$p(k,x)$ as the least integer greater or equal $(n-\frac{1}{k})x$. Then
$$R_{\sigma_i}*R_{\sigma_j}^{*k}*R_{\sigma_i}^{*x}*R_{\sigma_j}^{*p(k,x)}=R_{\sigma_i}^{*(x+1)}*R_{\sigma_j}^{*(p(k,x)+k)}.$$
\end{corollary}

Now we come back to the general analysis of the structure of ${\cal C}(Q)$.\\[2ex]
Recall that a dimension vector $d\in{\bf N}I$ is called a Schur root if there
exists a representation $X\in R_d$ such that ${\rm End}_{kQ}(X)\simeq k$. The
Schur roots play a prominent role in the study of generic properties of the
varieties $R_d$.\\[1ex]
Define the generic decomposition of a dimension vector
$d$ to be the unique (unordered) tuple $(d_1,\ldots,d_s)$ of dimension vectors
such that representations of the form
$$X_1\oplus\ldots\oplus X_s,\;\;\; X_k\in R_{d_k}\mbox{ for all }k$$
form a dense subset of $R_d$.\\[1ex]
By (\cite{Sc}, Theorem 2.1), a tuple $(d_1,\ldots,d_s)$ is the generic decomposition of $d$ if and only if:
\begin{enumerate}
\item $d=d_1+\ldots+d_s$,
\item each $d_k$ is a Schur root,
\item ${\rm ext}(d_k,d_l)=0$ for all $k\not=l$.
\end{enumerate}
This criterion has an immediate application to the structure of ${\cal C}(Q)$.

\begin{lemma}\label{cd} If $(d_1,\ldots,d_s)$ is the generic decomposition of $d$, then
$R_d$ is a commuting product of the $R_{d_k}$.
\end{lemma}

\proof We only have to note that, by Proposition \ref{obs}, we have
$$R_{d_k}*R_{d_l}=R_{d_k+d_l}=R_{d_l}*R_{d_k}$$
since all ${\rm ext}(d_k,d_l)$ vanish. \hb

\question Is the canonical decomposition already
characterized by this algebraic property?\\[3ex]
As a more important application of the notion of a Schur root to the structure
of ${\cal C}(Q)$, we will now derive a 'partial normal form' for its elements.

\begin{theorem}\label{pnf} Each element ${\cal A}\in{\cal C}(Q)$ can be written in the form
$${\cal A}=R_{d_1}*\ldots*R_{d_s},$$
where each $d_k$ is a Schur root, and
$${\rm ext}(d_k,d_{k+1})\not=0\mbox{ implies }{\rm ext}(d_{k+1},d_k)\not=0$$
for all $k=1\ldots s-1$.
\end{theorem}

\proof By definition, the element ${\cal A}$ can be written as
$${\cal A}=R_{d_1}*\ldots*R_{d_s}$$
for a tuple of Schur roots $d_*=(d_1,\ldots,d_s)$ (for example, we can take
$d_k=\sigma_{i_k}$ for some $i_k\in I$). We proceed by induction over
$$N(d_*):=\sum_{k<l}{\rm ext}(d_k,d_l).$$
If $N(d_*)=0$, there is nothing to prove. So, assume the theorem is
proved for all tuples $d_*$ such that $N(d_*)<N$ for a fixed $N\in{\bf N}$.\\[1ex]
If possible, we choose a position $k$ such that
$${\rm ext}(d_k,d_{k+1})\not=0\mbox{, but }{\rm ext}(d_{k+1},d_k)=0.$$
If this is not possible, the desired normal
form is already reached and we are done. Let $f_1,\ldots,f_t$ be the generic
decomposition of the dimension vector $d_k+d_{k+1}$. By Proposition \ref{obs}
and Lemma \ref{cd}, we have:
$$R_{d_k}*R_{d_{k+1}}=R_{d_k+d_{k+1}}=R_{f_1}*\ldots*R_{f_t}.$$
Thus, we can replace the tuple $d_*$ by
$$\widetilde{d}_*:=(d_1,\ldots,d_{k-1},f_1,\ldots,f_t,d_{k+2},\ldots,d_s)$$
without altering ${\cal A}$.
If we can prove that $N(\widetilde{d}_*)<N(d_*)$, we are done by
induction. To do this, we calculate:
\begin{eqnarray*} N(\widetilde{d}_*)-N(d_*)&=&
\sum_{p=1}^{k-1}\sum_{q=1}^t({\rm ext}(d_p,f_q)-{\rm ext}(d_p,d_k)-{\rm ext}(d_p,d_{k+1}))+\\
&&+\sum_{p=k+2}^{s}\sum_{q=1}^t({\rm ext}(f_q,d_p)-{\rm ext}(d_k,d_p)-{\rm ext}(d_{k+1},d_p))\\
&&-{\rm ext}(d_k,d_{k+1}).
\end{eqnarray*}
To estimate the summands on the right hand side, we first apply Proposition \ref{extest} to ${\cal A}=R_{d_p}$, ${\cal B}=R_{d_k}$, ${\cal C}=R_{d_{k+1}}$
for $p=1\ldots k-1$. This yields
$${\rm ext}(d_p,d_k+d_{k+1})\leq {\rm ext}(d_p,d_k)+{\rm ext}(d_p,d_{k+1}).$$
Next, we apply the same proposition repeatedly to ${\cal A}=R_{d_p}$, ${\cal B}=R_{f_q}$, ${\cal C}=R_{f_{q'}}$. Since all ${\rm ext}(f_q,f_{q'})$ for $q\not=q'$ vanish (by the properties of the generic decomposition), we get
the equality
$$\sum_{q=1}^t{\rm ext}(d_p,f_q)={\rm ext}(d_p,d_k+d_{k+1})\geq {\rm ext}(d_p,d_k)+{\rm ext}(d_p,d_{k+1}).$$
The same works for the second summand by duality. Moreover, ${\rm ext}(d_k,d_{k+1})$ is non-zero by assumption. The desired inequality
$N(\widetilde{d}_*)<N(d_*)$ is proved and we are done. \hb

\remark The above theorem gives only a 'partial normal form' in at least
two respects: first, in case ${\rm ext}(d_k,d_{k+1})=0={\rm ext}(d_k,d_{k+1})$,
the elements $R_{d_k}$ and $R_{d_{k+1}}$ commute, and there is no obvious choice
for their ordering in a normal form. Second, and more important, we only use
the relations provided by Proposition \ref{obs} to obtain this normal form.
But it is by no means clear in which cases these relations suffice to present
${\cal C}(Q)$.\\[1ex] Nevertheless, everything can be made to work in the
Dynkin case: in this case, we can choose an ordering on the finite set of positive
roots $d_1,\ldots,d_s$ (which are automatically Schur roots) such that
${\rm ext}(d_k,k_l)$ vanishes for $k\leq l$ (since the category of representations is directed). Using this notation, ${\cal C}(Q)$ coincides with the set of products
$$R_{d_1}^{*a_1}*\ldots*R_{d_s}^{*a_s}$$
for arbitrary $a_1,\ldots,a_s\in{\bf N}$. This follows easily from Theorem \ref{rqgd}.\\[2ex]
Completing our present discussion of the structure of ${\cal C}(Q)$, we end this section with the following questions:\\[2ex]
\questions $\;$
\begin{enumerate}
\item For which non-Dynkin quivers is the monoid ${\cal C}(Q)$ already presented
by the relations provided by Proposition \ref{obs}?
\item Which of these relations are 'essential' for general $Q$?
\item Is there an analogue of the 'partial normal form' of Theorem \ref{pnf}
in the degenerate quantized enveloping algebra ${\cal U}^+_0(Q)$?
\item If so, does this provide a tool for measuring the failure of $\eta$ to
be an isomorphism and the sufficiency of the relations \ref{obs}?
\end{enumerate}

\end{document}